\DeclareFontFamily{U}{mathb}{\hyphenchar\font45}
\DeclareFontShape{U}{mathb}{m}{n}{ <-6> matha5 <6-7> matha6 <7-8>
mathb7 <8-9> mathb8 <9-10> mathb9 <10-12> mathb10 <12-> mathb12 }{}
\DeclareSymbolFont{mathb}{U}{mathb}{m}{n}

\DeclareMathAccent{\abxring}{0}{mathb}{"38}

\DeclareFontFamily{U}{mathb}{\hyphenchar\font45}
\DeclareFontShape{U}{mathb}{m}{n}{ <-6> matha5 <6-7> matha6 <7-8>
mathb7 <8-9> mathb8 <9-10> mathb9 <10-12> mathb10 <12-> mathb12 }{}
\DeclareSymbolFont{mathb}{U}{mathb}{m}{n}

\documentclass[graybox]{svmult}

\usepackage{type1cm}        
\usepackage{makeidx}         
\usepackage{graphicx}        
\usepackage{multicol}        
\usepackage[bottom]{footmisc}

\usepackage{newtxtext}       %
\usepackage{newtxmath}       

\makeindex             

\usepackage{wrapfig}
\usepackage{graphicx}  
\usepackage{dcolumn}  
\usepackage{bm}       
\usepackage{algorithm,algpseudocode}
\usepackage{enumerate}
 \usepackage{booktabs}
\usepackage{multirow}
\usepackage{makecell}
\usepackage{pgfplots}
\usepackage{tikz}
\usetikzlibrary{patterns}
\usepackage{empheq}
\usepackage{todonotes}
\usepackage{xargs}
\usetikzlibrary{arrows} 
\usepackage{calrsfs}
\usepackage{enumitem}
\usepackage{upgreek}

\definecolor{myblack}{RGB}{53, 53, 53}
\definecolor{myblue}{RGB}{40, 75, 200}
\definecolor{myred}{RGB}{192, 50, 33}
\definecolor{myyellow}{RGB}{255, 166, 48}
\definecolor{mywhite}{RGB}{240, 237, 238}
\definecolor{mygreen}{RGB}{0, 102, 0}
\definecolor{mypurple}{RGB}{150, 0, 180}

\definecolor{green1}{RGB}{9, 82, 86}
\definecolor{green2}{RGB}{8, 127, 140}
\definecolor{green3}{RGB}{6, 167, 125}
\definecolor{green4}{RGB}{79, 109, 122}
\definecolor{green5}{RGB}{192, 214, 223}
\definecolor{violet}{RGB}{26,69,131}

\def \0{\mathbf{0}}			%

\def \x{{\mathbf{x}}}			
\def \s{{\mathbf{p}}}			
\def \R{{\mathbb{R}}}			
\def \N {\mathbb{N}}			
\def \vlb{{\mathbf{vl}}}			%
\def \vub{{\mathbf{vu}}}			%
\def \lb{{\mathbf{l}}}			%
\def \ub{{\mathbf{u}}}			%
\def \tu{{\mathbf{tu}}}			%
\def \tl{{\mathbf{tl}}}			%

\def \l{{\mathbf{l}}}			
\def \u{{\mathbf{u}}}			
\def \h{{h}}

\def \Rev{{\mathbf{R}}}			
\def \P{{\mathbf{P}}}			
\def \I{{\mathbf{I}}}			

\def \x{\vec{x}}				
\def \s{\vec{s}}				
\def \0{\vec{0}}			%

\def \R{\mathbb{R}}			
\def \N {\mathbb{N}}			

\def \Rev{\mathbf{R}}			
\def \P{\mathbf{P}}			
\def \I{\mathbf{I}}			

\DeclareMathAlphabet{\mathbfsf}{\encodingdefault}{\sfdefault}{bx}{n}
\DeclareMathAlphabet{\mathsfbfit}{\encodingdefault}{\sfdefault}{bx}{sl}
\DeclareMathAlphabet{\mathbfit}{\encodingdefault}{\rmdefault}{b}{it}
\DeclareMathAlphabet{\mathsfit}{\encodingdefault}{\sfdefault}{m}{sl}

\renewcommand{\vec}[1]{\boldsymbol{#1}}

\usepackage{tikz}
\usetikzlibrary{automata, positioning, arrows, external, matrix}
\usepackage{pgfplots}
\usepgfplotslibrary{fillbetween}
\usepgfplotslibrary{groupplots}
\usepgfplotslibrary{external}

\DeclareMathAlphabet{\pazocal}{OMS}{zplm}{m}{n}
\DeclareMathAlphabet{\bpazocal}{OMS}{cmsy}{b}{n}

\DeclareFontFamily{OT1}{pzc}{}
\DeclareFontShape{OT1}{pzc}{m}{it}{<-> s * [1.10] pzcmi7t}{}
\DeclareMathAlphabet{\mathpzc}{OT1}{pzc}{m}{it}

\def \F{\pazocal{F}} 
\def \A{\pazocal{A}} 
\def \L{\pazocal{L}} 
\def \S{\pazocal{S}} 
\def \X{\pazocal{X}} 

\algnewcommand\algorithmiconput{\textbf{Constants:}}
\algnewcommand\algorithmicinput{\textbf{Input:}}
\algnewcommand\algorithmicoutput{\textbf{Output:}}
\algnewcommand{\algorithmicgoto}{\textbf{go to}}%

\algnewcommand\Constants{\item[\algorithmiconput]}
\algnewcommand\Input{\item[\algorithmicinput]}%
\algnewcommand\Output{\item[\algorithmicoutput]}%
\algnewcommand{\Goto}[1]{\algorithmicgoto~\ref{#1}}%

\algrenewcommand\ALG@beginalgorithmic{\footnotesize}
\algrenewcommand\alglinenumber[1]{\scriptsize #1:}

\usetikzlibrary{decorations.pathmorphing,shapes}

\begin{document}

\title*{
Multilevel Active-Set Trust-Region (MASTR)  Method for Bound Constrained Minimization
}
\titlerunning{Multilevel Active-Set Trust-Region Method} 
\author{Alena Kopani\v{c}\'akov\'a and Rolf Krause}
\institute{Alena Kopani\v{c}\'akov\'a \at Universit\`{a} della Svizzera italiana, Switzerland,  \email{alena.kopanicakova@usi.ch}
\and Rolf Krause \at Universit\`{a} della Svizzera italiana, Switzerland,   \email{rolf.krause@usi.ch}}
%
%
\maketitle

\abstract*{We introduce a novel variant of the recursive multilevel trust-region (RMTR) method, called MASTR. 
The method is designed for solving non-convex bound-constrained minimization problems, which arise from the finite element discretization of partial differential equations.
MASTR utilizes an active-set strategy based on the truncated basis approach in order to preserve the variable bounds defined on the finest level by the coarser levels. 
Usage of this approach allows for fast convergence of the MASTR method, especially once the exact active-set is detected. 
The efficiency of the method is demonstrated by means of numerical examples.  }

\abstract{We introduce a novel variant of the recursive multilevel trust-region (RMTR) method, called MASTR. 
The method is designed for solving non-convex bound-constrained minimization problems, which arise from the finite element discretization of partial differential equations.
MASTR utilizes an active-set strategy based on the truncated basis approach in order to preserve the variable bounds defined on the finest level by the coarser levels. 
Usage of this approach allows for fast convergence of the MASTR method, especially once the exact active-set is detected. 
The efficiency of the method is demonstrated by means of numerical examples. }

\section{Introduction}
\label{sec:intro}
We consider a minimization problem of the following type:
\begin{equation}
\begin{aligned}
\underset{\x \in \R^n}{\text{min}} \quad  &  f(\x) \\
\text{subject to} \quad &\x \in \pazocal{F},
\end{aligned}
\label{eq:problem_min_tr_bla}
\tag{P}
\end{equation}
where $f: \R^n \rightarrow \R$ is possibly non-convex but continuously differentiable.
The feasible set $\pazocal{F}:=\{\x \in \R^{n} | \ \l \leq \x \leq \u \}$ is defined in terms of the pointwise lower bound $\lb \in \R^n$ and the upper bound $\ub \in \R^n$.
We assume that the function $f$ arises from the finite element (FEM) discretization of partial differential equations (PDEs). 
Here, $n \in \N$ denotes the dimension of the finite element space and it is typically very large. 
Problems of this type arise commonly in many scientific applications, for example in fracture or contact mechanics~\cite{kopanivcakova2020recursive, krause2009nonsmooth}.  

Multilevel methods are known to be optimal solution strategies for systems arising from the discretization of, usually elliptic, PDEs, as their convergence rate is often independent of the problem size and the number of required arithmetic operations grows proportionally with the number of unknowns.
These methods have been originally designed for unconstrained PDEs~\cite{briggs2000multigrid}. 
Their extension to constrained settings is not trivial as the coarse levels are often not capable of resolving the finest-level constraints sufficiently well, especially if the constraints are oscillatory~\cite{kornhuber2001adaptive}. 
The initial attempts to incorporate the constraints into the multilevel framework are associated with solving linear complementarity problems, see for instance~\cite{mandel1984multilevel, brandt1983multigrid, hackbusch1983multi, gelman1990multilevel}. 
The devised methods employed various constraint projection rules for constructing the coarse-level variable bounds, such that coarse-level corrections are admissible by the finest level. 
Unfortunately, these projection rules provided quite a narrow approximation of the finest-level constraints.
As a consequence, the resulting multilevel methods converge significantly slower than standard linear multigrid.
In order to enhance the convergence speed, Kornhuber proposed an active-set multigrid method~\cite{kornhuber1994monotone}. 
The method utilizes a truncated basis approach and recovers the convergence rate of the unconstrained multigrid, once the exact active-set is detected~\cite{hoppe1994adaptive, kornhuber1994monotone}. 

In the field of nonlinear optimization, a very few existing nonlinear multilevel algorithms can be readily employed. 
For instance, Vallejos proposed a gradient projection based multilevel method~\cite{vallejos2010mgopt}. 
Two multilevel line-search methods, designed for convex optimization problems, are proposed in~\cite{kovcvara2016first}.
These methods utilize constraint projection rules developed in~\cite{hackbusch1983multi} and a variant of the active-set strategy from~\cite{kornhuber1994monotone}. 
In the context of non-convex optimization problems, Gratton et~al. proposed a variant of the recursive multilevel trust-region (RMTR) method~\cite{gratton2008_inf} by utilizing the constraint projection rules from~\cite{gelman1990multilevel}. 
To our knowledge, this is currently the only nonlinear multilevel method, which provides global convergence guarantees for non-convex bound constrained optimization problems. 

In the presented work, we propose to enhance the convergence speed of the RMTR method~\cite{gratton2008_inf}.
More precisely, we introduce an active-set variant, called MASTR, which utilizes the  truncated basis approach~\cite{kornhuber1994monotone}. 
In contrast to~\cite{kovcvara2016first}, we employ coarse-level models of the Galerkin type.
This simplifies the practical implementation of the algorithm and avoids unconventional modifications to existing FEM software packages. 
As it will be demonstrated by our numerical results, employing the active-set approach allows for significant speedup of the RMTR method. 

\section{Recursive multilevel trust-region (RMTR) method}
\label{sec:RMTR}
In this work, we minimize~\eqref{eq:problem_min_tr_bla} using a novel variant of the RMTR method~\cite{gratton2008_inf}. 
RMTR combines the global convergence properties of the trust-region (TR) method with the efficiency of multilevel methods.
By design, the RMTR method employs a hierarchy of $L$ levels. 
Each level~$l$ is associated with mesh~$\pazocal{T}^l$, which encapsulates the computational domain $\Omega \in \R^d$, where $d\in \N$.
The mesh~$\pazocal{T}^l$ is used to construct the finite-dimensional FEM space~$\X^l$, spanned by the basis functions~$\{ N_k^l \}_{k \in \pazocal{N}^l}$, where~$\pazocal{N}^l$ denotes the set of interior nodes of the mesh~$\pazocal{T}^l$.  
The support of a given basis function $N_k^l$ is defined as $\omega^l_k = \overline{ \{ x \in \Omega \ | \ N_k^l(x) \neq 0 \}}$. 
Note, the support is typically local and restricted only to the neighborhood of the~{$k$-th node} of the mesh~$\pazocal{T}^l$. 

The transfer of data between subsequent levels of the multilevel hierarchy is carried out using three transfer operators, namely 
prolongation $\I_l^{l+1}: \R^{n^l} \rightarrow \R^{n^{l+1}}$, restriction $\Rev_{l+1}^l:= (\I_l^{l+1})^T $ and projection $\P_{l+1}^l: \R^{n^{l+1}} \rightarrow \R^{n^l}$. 

\subsection{Algorithm}
On each level $l$, the RMTR method approximates~\eqref{eq:problem_min_tr_bla} by means of some level-dependent objective function $h^l:\R^{n^l} \rightarrow \R$ and feasible set $\pazocal{F}^l:=\{\x^l \in \R^{n^l} | \ \l^l \leq \x^l \leq \u^l \}$.
The function $h^l$ is approximately minimized in order to obtain coarse-level correction. 
This correction is then interpolated to the subsequent finer level, where it is used to improve the current iterate.

More precisely, the algorithm starts on the finest level, $l=L$, with an initial iterate $\x_0^L$ and passes through all levels until the coarsest level, $l=1$, is reached. 
On each level $l$, the algorithm performs $\mu_1$ pre-smoothing steps to improve the current iterate $\x_{0}^l$. 
The smoothing is performed using the TR method~\cite{conn2000trust}.
Thus, on each TR iteration $i$, the search direction $\s_i^l$ is obtained by approximately solving following minimization problem:
\begin{align}
\underset{\s_i^l \in \R^{n^l}}{\text{min}} \ m_i^l(\s_i^l ) := &h^l(\x_i^l) + \langle \nabla h^l(\x_i^l), \s_i^l \rangle + \frac{1}{2} \langle \s_i^l,  \nabla^2 h^l(\x^l_i) \  \s^l_i \rangle, \nonumber  \\
&\text{such that}  \ \ \x^l_i + \s^l_i \in  \pazocal{F}^l, \label{eq:model_qp}  \\
& \| \s_i^l \|_{\infty} \leq \Delta_i^l, \nonumber
\end{align}
where $m_i$ is second-order Taylor approximation of $h^l$. 
The symbol $\Delta_i^l > 0$ denotes a TR radius, which controls the size of the correction $\s_i^l$. 
In contrast to line-search methods, the correction $\s_i^l$ is used only, if $\rho_i^l > \eta_1$, where $ \rho_i^l =  \frac{h^l(\x_i) - h^l(\x_i + \s^l_i)}{m^l(\s^l_i)}$ and $\eta_1>0$. 
Otherwise, $\s^l_i$ is disposed and the size of the TR radius is reduced. 
The result of the pre-smoothing, the iterate $\x_{\mu_1}^l$,  is then used to initialize the solution vector on the subsequent coarser level, i.e.,~$\x_{0}^{l-1} := \P^{l-1}_l \x_{\mu_1}^l$.

Once the coarsest level is reached, we apply $\mu^1$ steps of the TR method to obtain the updated iterate $\x_{\mu^1}^1$. 
The algorithm then returns to the finest level. 
To this aim, the correction obtained on the level $l$, i.e.,~$\x_{\mu^{l}}^{l} - \x_{0}^{l}$, is transfered to the level~$l+1$. 
Here, the symbol $\mu^{l}$ denotes a sum of all iterations taken on a given level $l$.
However, the quality of the prolongated coarse-level correction $\s_{\mu_1+1}^{l+1} := \I^{l+1}_{l}(\x_{\mu^{l}}^{l} - \x_{0}^{l})$ has to be assessed before it is accepted on the level $l+1$.  
For this reason, we define a multilevel TR ratio as
\begin{align}
\rho^{l+1}_{\mu_1+1} := \frac{h^{l+1}(\x_{\mu_1}^{l+1}) - h^{l+1}(\x_{\mu_1}^{l+1} + \s_{\mu_1+1}^{l+1})}{h^{l}(\x_{0}^{l}) - h^{l}(\x_{\mu^{l}}^{l})}.
\label{eq:rho_ml}
\end{align}
The correction $\s_{\mu_1+1}^{l+1}$ is accepted if $\rho^{l+1}_{\mu_1+1} > \eta_1$. 
If $\rho^{l+1}_{\mu_1+1} \leq \eta_1$, the correction~$\s_{\mu_1+1}^{l+1}$ is rejected.
Additionally, the TR radius has to be updated accordingly. 
To this end, the RMTR algorithm performs $\mu_2$ post-smoothing steps at a given level $l$.
This process is repeated on every level until the finest level is reached.

\textcolor{white}{a}\\
\noindent \textbf{Construction of level-dependent objective functions and feasible sets}\\
In this work, we create a level-dependent objective function $\h^l$ as follows:
\begin{align}
\h^{l}(\x^{l}) :=  \langle \Rev^{l}_{l+1} \nabla h_{\mu_1}^{l+1},\x^l - \x^l_0 \rangle + \frac{1}{2} \langle \x^l - \x^l_0, ( \Rev^l_{l+1} \nabla^2 h_{ \mu_1}^{l+1} \I_l^{l+1})  (\x^l - \x^l_0) \rangle,
\label{eq:coarse_objective_galerkin}
\end{align}
where~$\Rev^{l}_{l+1} \nabla h_{\mu_1}^{l+1}$ and~$\Rev^l_{l+1} \nabla^2 h_{\mu_1}^{l+1} \I_l^{l+1}$ represent  the restricted gradient and the Hessian from the level~$l+1$, respectively. 
As we will see in Section~\ref{sec:active_set}, employing the coarse-level models of this particular type allows for straightforward incorporation of the active set strategy within the multilevel settings. 

The level-dependent feasible set~$\F^l$ is created by intersecting the set $\L^l$ with the set  $\S^l$, thus  as $\F^l := \L^l \cap \S^l$. 
The role of the set $\pazocal{S}^l := \{ \x^l \in \R^{n^l} \ | \ \tl^l \leq \x^l   \leq   \tu^l  \}$ is to ensure that the size of the prolongated coarse-level correction remains bounded by the TR radius~$\Delta^l_{\mu_1}$, i.e.,~$\| \I_{l}^{l+1} \s^l \|_{\infty} \leq \Delta^l_{\mu_1}$. 
To this aim, we construct $\pazocal{S}^l$ by employing the projection rules especially designed for TR bounds in~\cite{gratton2008_inf}.

The function of the set $\L^l := \{ \x^l \in \R^{n^l} \ | \ \vlb^l \leq \x^l \leq \vub^l \}$ is to guarantee that the  prolongated  coarse-level correction produces a feasible trial point, i.e.,~$\x_{\mu_1}^{l+1} +  \I_{l}^{l+1} \s^l  \in \pazocal{F}^{l+1}$. 
Following~\cite{gelman1990multilevel, gratton2008_inf}, we  can construct $\vlb^l, \vub^l$ in a component-wise manner as
\begin{equation}
\begin{aligned}
(\vlb^l)_k &:= (\x_{0}^l)_k + 
  \underset{j \in \pazocal{N}^{l+1} \cap \  \abxring{\omega}_k^l}{\max} [ (\vlb^{l+1} - \x^{l+1}_{\mu_1})_j ], \\ 
(\vub^l)_k &:= (\x_{0}^l)_k +
  \underset{j \in \pazocal{N}^{l+1} \cap \  \abxring{\omega}_k^l }{\min} [ (\vub^{l+1} - \x^{l+1}_{\mu_1})_j ],
\label{eq:hard_constraints_KV}
\end{aligned}
\end{equation}
where $(\cdot)_k$ denotes the $k$-th component of a given vector. 
Note, the support ${\omega_k^l}$ of the basis function $N_k^l$ (associated with~$k$-th node of the mesh~$\pazocal{T}^l$) determines, which components of the variable bounds~$\vlb^{l+1}$ and~$\vub^{l+1}$ have to be taken into account while constructing $\vlb^l, \vub^l$.

\begin{remark}
Throughout this work, we assume that $h^L:=f$ and $\pazocal{F}^L:=\pazocal{F}$. 
\end{remark}

\section{Multilevel active-set trust-region (MASTR) method}
\label{sec:active_set}
In this section, we present how to incorporate the active-set strategy into the RMTR framework.
The devised algorithm has also a form of the standard V-cycle. 
The key idea behind the proposed MASTR method is to identify an active-set
\begin{align}
\A^l :=\{ k \in \{1, \ldots, n^l \} \ | \ (\vlb^l)_k = (\x^l_{\mu_1})_k \ \text{or} \ (\vub^l)_k = (\x^l_{\mu_1})_k \},
\label{eq:active_set}
\end{align}
before descending to the coarser level. 
Here, the vectors $\vlb^l, \vub^l$ denote lower and upper bounds that define the set~$\L^l$, c.f.~\eqref{eq:hard_constraints_KV}.
The components of the solution vector~$\x^l_{\mu_1}$, which belong to the active-set $\A^{l}$, are then held fixed and cannot be altered by the coarser levels.
To this aim, the level-dependent objective functions $\{ \h^a \}_{a=1}^{l-1}$ and feasible sets $\{ \pazocal{L}^a \}_{a=1}^{l-1}$ have to be constructed such that the minimization process on a given level yields coarse-level corrections, which fulfil this requirement. 
Following~\cite{kornhuber1994monotone, hoppe1994adaptive}, we construct $\{ \h^a \}_{a=1}^{l-1}$ and $\{ \pazocal{L}^a \}_{a=1}^{l-1}$ using a truncated basis method.

\textcolor{white}{a}\\
\noindent \textbf{Construction of truncated FEM spaces}\\
The truncated basis method~\cite{kornhuber1994monotone} constructs truncated FEM spaces~$\{ \widetilde{\X}^l \}_{l=1}^{L-1}$ by exploiting the fact that the 
basis functions on level~$l$, can be written as a linear combination of basis functions on level~$l+1$, i.e.~$ N_k^l = \sum^{n^{l+1}}_{p=1} (\I_{l}^{l+1})_{pk} N_p^{l+1}$.
Note, this property is also utilized while constructing coarse level model of Galerkin~type,~c.f.~\eqref{eq:coarse_objective_galerkin}. 

Each truncated FEM space~$\widetilde{\X}^l$ is spanned by truncated basis functions~$\{ \widetilde{N}_k^l \}_{k \in \pazocal{N}^l}$, the value of which is zero at all active nodes of the finer levels. 
More precisely, we can construct truncated basis functions in a recursive manner as
\begin{align}
\widetilde{N}_k^{l} = \sum^{n^{l+1}}_{p=1} (\widetilde{\I}_{l}^{l+1})_{pk}  \widetilde{N}_p^{l+1},
 \label{eq:basis_recursion_truncated_mat}
\end{align}
where~$\widetilde{\I}_{l}^{l+1}$ is truncated prolongation operator defined by
\begin{align}
(\widetilde{\I}_{l}^{l+1})_{pk} = 
\begin{cases}
0,  \quad &\text{if }\ p \in \A^{l+1}, \\
({\I}_{l}^{l+1})_{pk}, \quad& \text{otherwise}.
\end{cases}
\label{eq:truncation}
\end{align}
The operator~$\widetilde{\I}_{l}^{l+1}$ is obtained from the prolongation operator~${\I}_{l}^{l+1}$ by setting~$p$-th row of~${\I}_{l}^{l+1}$ to zero, for all~$p \in \A^{l+1}$.
The application of~$\widetilde{\I}_{l}^{l+1}$ in~\eqref{eq:basis_recursion_truncated_mat} removes contributions of basis functions associated with active nodes on level $l+1$, defined by the set~$\A^{l+1}$.

\begin{remark}
The functions $\{ \widetilde{N}_k^{L-1} \}_{k \in \pazocal{N}^{L-1}}$ are constructed using~\eqref{eq:basis_recursion_truncated_mat} with $\{ {N}_k^{L} \}_{k \in \pazocal{N}^{L}}$.
\end{remark}

\textcolor{white}{a}\\
\noindent \textbf{Construction of level-dependent objective functions and feasible sets}\\
Using  truncated FEM spaces~$\{ \widetilde{\X}^l \}_{l=1}^{L-1}$, we can now construct level-dependent objective functions $\{ \h^l \}_{l=1}^{L-1}$ and feasible sets $\{ \pazocal{F}^l\}_{l=1}^{L-1}$.
In particular, for a given level $l<L$,  the level-dependent objective function $\h^{l}: \R^{n^{l}} \rightarrow \R$ is created as follows:
\begin{align}
\h^{l}(\x^{l}) :=  \langle(\widetilde{\I}_{l}^{l+1})^T  \nabla h_{\mu_1}^{l+1}, \x^l - \x^l_0 \rangle + \frac{1}{2} \langle   \x^l - \x^l_0, ( \widetilde{\I}_{l}^{l+1})^T \nabla^2 h_{ \mu_1}^{l+1} \widetilde{\I}_{l}^{l+1})  (\x^l - \x^l_0) \rangle,
\label{eq:galerkin_trunc}
\end{align}
where we used truncated transfer operator~$\widetilde{\I}_{l}^{l+1}$ to restrict gradient~$ \nabla h_{\mu_1}^{l+1}$ and Hessian~$\nabla^2 h_{ \mu_1}^{l+1}$ from level $l+1$ to level $l$.
The application of~$\widetilde{\I}_{l}^{l+1}$ in~\eqref{eq:galerkin_trunc} removes the components of fine-level gradient/Hessian associated with the active-set $\A^{l+1}$.
Please note, the formulation~\eqref{eq:galerkin_trunc}  does not require explicit representation of~$\{ \widetilde{\X}^l \}_{l=1}^{L-1}$. 

The construction of each level-dependent feasible set~$\L^l$ can be performed using projection rules defined by~\eqref{eq:hard_constraints_KV}. 
However, formulas~\eqref{eq:hard_constraints_KV} are now determined by the support of the truncated basis functions, spanning~$\widetilde{\X}^l$. 
Since the support of basis functions spanning~$\widetilde{\X}^l$ is different from the support of the basis functions spanning~${\X}^l$, fewer components of a fine-level variable bounds are taken into account by~\eqref{eq:hard_constraints_KV}. 
This yields less restrictive coarse-level constraints and allows for larger coarse grid corrections.
Algorithm~\ref{alg:rmtr_active_set} summarizes the proposed MASTR method. 

\begin{algorithm}[t]
\footnotesize
\caption{\footnotesize MASTR($l, \  h^l, \  \x^l_{0}, \ \F^l, \  \Delta^l_{0}$)}
\label{alg:rmtr_active_set}
\begin{algorithmic}[1]
\Require{$ l \in \N, \  h^l:\R^{n^l}\rightarrow \R,  \ \x^l_{0} \in \R^{n^l}, \ \F^l, \  \Delta^l_{0} \in \R$}
\Constants { $\mu_1, \ \mu_2, \ \mu^1 \in \N$}
\State $ [\x_{\mu_1}^l, \  \Delta_{\mu_1}^l] = \text{Trust\_region}(\h^l,  \ \x_{0}^l, \  \F^l,  \  \Delta_{0}^l, \  \mu_1 )$
\Comment{Pre-smoothing}
\State Construct~$\A^l$, $\widetilde{\I}_{l-1}^l$, $\h^{l-1}$,  $\F^{l-1}$ \Comment{\small Initialize coarse-level quantities}
\If{$l ==2$}
\State $ [\x_{\mu^{l-1}}^l,  \_\_] = \text{Trust\_region}(\h^{l-1},  \ \P_{l}^{l-1} \x_{\mu_1}^l, \  \F^{l-1},  \  \Delta_{\mu_1}^{l}, \  \mu^1 )$  \Comment{Coarse-level solve}
\Else
\State $[\x_{\mu^{l-1}}^{l-1}, \_\_]$=MASTR($l-1, \h^{l-1},  \P_{l}^{l-1} \x_{\mu_1}^l, \  \F^{l-1}, \ \Delta_{\mu_1}^l) $  \Comment{Call MASTR recursively}
\EndIf
\State $\s^l_{\mu_1+1} = \widetilde{\I}_{l-1}^{l}(\x_{\mu^{l-1}}^{l-1} - \P_{l}^{l-1} \x_{\mu_1}^l)$  \Comment{Prolongate coarse-level correction}
\State Compute $\rho^l_{\mu_1 +1}$ by means of~\eqref{eq:rho_ml}
\State $ [\x_{\mu_1+1}^l, \  \Delta_{\mu_1+1}^l]$ = Convergence\_control($\rho^l_{\mu_1 +1}, \ \x_{\mu_1}^l, \ \s^l_{\mu_1+1}, \ \Delta_{\mu_1}^l$) 
 \Comment{\footnotesize Update iterate and tr. radius}
\State $ [\x_{*}^l, \  \Delta_{*}^l] = \text{Trust\_region}(\h^l,  \ \x_{\mu_1+1}^l, \  \F^l,  \  \Delta_{\mu_1+1}^l, \  \mu_2 )$
 \Comment{ \footnotesize Post-smoothing}
\State \Return $\x_{*}^l$, $\Delta_*^l$ 
\end{algorithmic}
\end{algorithm}

\section{Numerical results}
\label{sec:num_results}
We study the performance of the proposed MASTR method using three numerical examples. 
Examples are defined on domain~${\Omega:=[0,1]^2}$ with boundary~${\Gamma=\partial \Omega}$, decomposed into three parts:~${\Gamma_l = \{0\} \times [0,1]}$,~${\Gamma_r = \{1\} \times [0,1]}$, and~${\Gamma_f = [0,1] \times \{0,1\}}$.  
The discretization is performed using uniform mesh and~$\mathbb{Q}_1$ Lagrange finite elements.

\noindent\textbf{Ex.1. MEMBRANE:}
Let us consider the following  minimization problem~\cite{domoradova2007projector}: 
\begin{equation}
\begin{aligned}
& \underset{u \in \X}{\text{min}} \ \ f(u) := \frac{1}{2} \int_{\Omega} \| \nabla u(x) \|^2 \ d x + \int_{\Omega} u(x) \ d x, \\
& \text{subject to} \ \ \text{lb}(x) \leq u, \quad \text{on~$\Gamma_r$}.
\end{aligned}
\label{eq:min_membrane}
\end{equation}
The lower bound~$\text{lb}$ is defined on the right part of the boundary,~$\Gamma_r$, by the upper part of the circle with the radius,~$r=1$, and the center,~$C=(1;-0.5;-1.3)$.
The minimization in~\eqref{eq:min_membrane} is performed over the space~$\X:=\{u \in H^1(\Omega) \ | \ u = 0 \ \text{on} \ \Gamma_l \}$.

\noindent\textbf{Ex.2. IGNITION:}
Following~\cite{briggs2000multigrid, kovcvara2016first}, we minimize following optimization problem:
\begin{equation}
\begin{aligned}
& \underset{u \in \X}{\text{min}} \ \  f(u) := \frac{1}{2} \int_{\Omega} \| \nabla u(x) \|^2 - (u e^u - e^u) \ d x - \int_{\Omega} f(x) u \ d x, \\
& \text{subject to} \ \  \text{lb}(x) \leq u \leq \text{ub}(x), \quad  \text{a.e. in~$\Omega$}.
\end{aligned}
\label{eq:min_combustion}
\end{equation}
The variable bounds and right-hand side are defined as
\begin{equation*}
\begin{aligned}
\text{lb}(x) &= -8(x_1 - 7/16)^2 - 8(x_2 - 7/16)^2 +0.2, \quad \text{ub}(x) = 0.5,  \\
f(x) &= (9 \pi^2 + e^{(x_1^2-x_1^3)\sin(3\pi x_2)} (x_1^2-x_1^3) + 6 x_1 - 2)\sin(3\pi x_1),
\end{aligned}
\end{equation*}
where~$f \in L^2(\Omega)$ and~$x_1, x_2$ denote spatial coordinates. 
The minimization~\eqref{eq:min_combustion} is carried out over the space~$\X:=\{u \in H^1(\Omega) \ | \ u = 0 \ \text{on} \ \Gamma \}$.

\noindent\textbf{Ex.3. MOREBV:}
We consider the following non-convex minimization problem~\cite{gratton2010numerical}:
\begin{equation}
\begin{aligned}
& \underset{u \in \X}{\text{min}} \ \  f(u) :=  \int_{\Omega}  \| \Delta u(x)  - 0.5(u(x) + \langle {e}, x \rangle + 1)^3  \|_2^2 \ d x, \\
& \text{subject to} \ \  \text{lb}(x) \leq u, \quad  \text{a.e. in~$\Omega$}, 
\end{aligned}
\label{eq:min_morevb}
\end{equation}
where ${e}$ denotes a unit vector. The lower bound is defined as
\begin{equation*}
\begin{aligned}
\text{lb}(x) = \sin(5 \pi x_1) \sin(\pi x_2) \sin(\pi (1-x_1)) \sin( \pi (1-x_2)).
\end{aligned}
\end{equation*}
where~$x_1, x_2$ denote spatial coordinates. 
The minimization in~\eqref{eq:min_morevb} is performed over the space~$\X:=\{u \in H^1(\Omega) \ | \ u = 0 \ \text{on} \ \Gamma \}$.

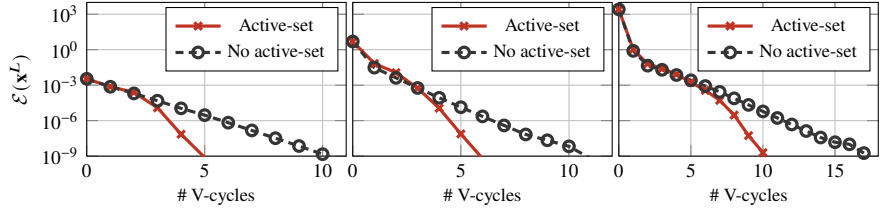
\begin{figure}[t]
\centering
 \begin{tikzpicture}[]
    \begin{groupplot}[
        group style={
            group size = 3 by 1,
            horizontal sep=2.5pt,
            vertical sep=3cm,
          },
	legend pos=north east,
	ytick={1e3, 1, 1e-3, 1e-6,  1e-9},
        width=0.43\textwidth, 
        height=0.31\textwidth, 
	grid=major, 
	xmode=normal,
	ymode=log,
	xlabel={\scriptsize \# V-cycles},
	x label style={at={(axis description cs:0.5, 0.1)}},	
	ymin=1e-9, 
	ymax=1e4,
	xmin=0,	
	tick label style={font=\scriptsize},
	label style={font=\scriptsize},
	legend style={font=\scriptsize}        
      ]

      \nextgroupplot[align=left, 
      title={}, 
      y label style={at={(axis description cs:0.1,.5)}},      
      ylabel={$\scriptsize \pazocal{E}(\x^L)$},   
      	legend pos=north east,
      legend entries={\scriptsize Active-set, \scriptsize No active-set},    
      xmax=11]
         \addplot[color = myred, very thick,  mark=x] table [x=its, y=g_trun, col sep=comma] {membrane.csv};
   \addplot[color = myblack, very thick, densely dashed, mark=o, mark options={solid}] table [x=its, y=g_loc, col sep=comma] {membrane.csv};

      \nextgroupplot[align=left, 
      title={}, 
      ylabel={},   
	yticklabels={},      
      	legend pos=north east,
      legend entries={\scriptsize Active-set, \scriptsize No active-set},    
      xmax=12]

   \addplot[color = myred, very thick,  mark=x] table [x=its, y=g_trun, col sep=comma] {ignition.csv};
   \addplot[color = myblack, very thick, densely dashed, mark=o, mark options={solid}] table [x=its, y=g_loc, col sep=comma] {ignition.csv};

      \nextgroupplot[align=left, 
      title={}, 
	yticklabels={},
	      	legend pos=north east,
      legend entries={\scriptsize Active-set, \scriptsize No active-set},    
      xmax=18]	

 \addplot[color = myred, very thick,  mark=x] table [x=its, y=g_trun, col sep=comma] {morebv.csv};
    \addplot[color = myblack, very thick, densely dashed, mark=o, mark options={solid}] table [x=its, y=g_loc, col sep=comma] {morebv.csv};      

    \end{groupplot}
  \end{tikzpicture}
  \caption[The convergence of the RMTR method with and without an active-set strategy.]{
The convergence of the RMTR method with (red color) and without (black color) an active-set strategy.  
 \emph{Left:} MEMBRANE.
 \emph{Middle:} IGNITION.
 \emph{Right:} MOREBV. 
  }
  \label{fig:test_constraints}
\end{figure}

\subsection{Convergence study}
We compare the convergence behavior of the the proposed MASTR method with the standard RMTR method (without the active-set strategy).
Both methods are implemented as part of the open-source library UTOPIA~\cite{utopia}. 
The performed study considers a setup with six~levels and one~pre/post-smoothing step. 
The trust-region subproblems~\eqref{eq:model_qp} are solved using one~iteration of successive coordinate minimization~\cite{gratton2008_inf}.
The algorithms terminate, if $\pazocal{E}(\x^L)<10^{-9}$ is satisfied.
The criticality measure~$\pazocal{E}(\x)$ is defined as $\pazocal{E}(\x):=\| \pazocal{P}(\x-\nabla f(\x)) - \x \|$, where $\pazocal{P}$ is the orthogonal projection onto the feasible set~$\pazocal{F}$.

As we can see from Figure~\ref{fig:test_constraints}, using an active-set approach is beneficial, as it allows for significant speed up.
We can also observe that during the active-set identification phase (first few V-cycles), both approaches are comparable. 
However,  once the exact active-set is detected, MASTR accelerates~and~converges~faster~than~standard~RMTR.


\end{document}